\documentclass[11pt]{article}
\usepackage{graphicx,mathrsfs,amssymb}
\usepackage{amsmath,latexsym,color, amsbsy,amssymb}
\usepackage{enumerate}
\usepackage{graphicx}
\usepackage{verbatim}

\usepackage[normalem]{ulem}
\usepackage{soul}

\makeatletter
\newcommand{\superimpose}[2]{%
	{\ooalign{$#1\@firstoftwo#2$\cr\hfil$#1\@secondoftwo#2$\hfil\cr}}}
\makeatother
\newcommand{\transversal}{\mkern-1mu\mathrel{\mathpalette\superimpose{{\top}{\scrunch{\cap}}}}\mkern-1mu}
\newcommand{\scrunch}[1]{\resizebox{\width}{.9\height}{$#1$}}

\def\ds{\displaystyle}
\def\forall{\hbox{for all}~}

\def\argmax{\hbox{arg}\!\max}

\def\ve{\varepsilon}
\def\n{\noindent}

\def\R{\mathbb{R}}

\def\vs{\vskip 2em}

\def\v{\vskip 1em}

\def\C{{\cal C}}
\def\bega{\begin{array}}
	\def\enda{\end{array}}
\def\begi{\begin{itemize}}
	\def\endi{\end{itemize}}

\def\I{{\mathcal I}}

\def\bel{\begin{equation}\label}
	\def\eeq{\end{equation}}
\def\sqr#1#2{\vbox{\hrule height .#2pt
		\hbox{\vrule width .#2pt height #1pt \kern #1pt
			\vrule width .#2pt}\hrule height .#2pt }}
\def\square{\sqr74}
\def\endproof{\hphantom{MM}\hfill\llap{$\square$}\goodbreak}

\newtheorem{theorem}{Theorem}[section]

\newtheorem{definition}[theorem]{Definition}
\newtheorem{remark}[theorem]{Remark}

\newtheorem{proposition}{Proposition}[section]

\begin{document}
	\title{\bf A lower bound on the quantitative version of the transversality theorem}\vs
	\author{Andrew Murdza and Khai T. Nguyen\\ 
		\\
		{\small Department of Mathematics, North Carolina State University}\\
		{\small e-mails: ~apmurdza@ncsu.edu,  ~ khai@math.ncsu.edu}
		%		\\
		%		\\
		%		{\small\it Dedicated to Professor Duong Minh Duc on the occasion of his 70th birthday}
	}
	\maketitle
	\begin{abstract} The present paper  studies a quantitative version of the transversality theorem. More precisely, given a continuous function $f\in \mathcal{C}([0,1]^d,\R^m)$ and a  manifold  $W\subset \R^m$ of dimension $p$,  a sharpness result on the upper quantitative estimate  of  the $(d+p-m)$-dimensional Hausdorff measure of the set $\mathcal{Z}_{W}^{f}=\left\{x\in [0,1]^d: f(x)\in W\right\}$, which was achieved in \cite{MN}, will be proved in terms of power functions. 
		\medskip
		
		\n {\bf AMS Mathematics Subject Classification.} 46T20
	\end{abstract}
	\section{Introduction}
	\label{sec:1}
	Let $g:X\to Y$ be a $\mathcal{C}^{1}$ map between two smooth manifolds $X$ of dimension $d$ and $Y$ of dimension $m$ . For any smooth submanifold $W\subseteq Y$ of dimension $p$, we say that the function $g$ is transverse to $W$ and write  $g \transversal W$ if
	\[
	(dg)_p(T_pX)+T_{g(p)}(W)~=~T_{g(p)}(Y)\quad\forall p\in g^{-1}(W).
	\]
	The transversality lemma, which is the key to studying  Thom's  transversality theorem \cite{Thom,N1,N2}, shows that the set of  transverse maps is dense \cite{Th1}. In particular,  for any continuous function $f:[0,1]^d\to\R^m$ and 
	any $\ve>0$, there exists a  $\mathcal{C}^{1}$  function $f_{\ve}: [0,1]^d\to\R^m$ such that 
	\[
	\left\|f_{\ve}-f\right\|_{\mathcal{C}^1}~\leq~\ve\qquad \mathrm{and}\qquad f_{\ve}\transversal W.
	\]
	For every $h\in\mathcal{C}\big([0,1]^d,\R^m\big)$, consider the set
	\bel{Z}
	\mathcal{Z}^{h}_{W}~:=~\left\{x\in [0,1]^d: h(x)\in W\right\}.
	\eeq
	If $h$ is smooth and transverse to $W$, 
	then $\mathcal{Z}^h_{W}$ is a  $(d+p-m)$-dimensional smooth manifold.  Hence, its $(d+p-m)$-dimensional 
	Hausdorff measure is finite.   In  the spirit of  metric entropy, which was used in the study of compactness estimates for solution sets of hyperbolic conservation laws  \cite{AGN1,AGN2,AGN,CDN} and Hamilton-Jacobi equations \cite{ACN1, ACN2,BDN}, a natural question is to  perform a quantitative analysis of the  measure  of $\mathcal{Z}^{f}_{W}$.  Namely, how small can one
	make this measure, by an $\ve$-perturbation of $f$? To formulate more precisely the result, given $f\in\mathcal{C}([0,1]^d,\R^m)$, one defines
	\bel{N-f}
	\mathcal{N}^{f}_{W}(\ve)~:=~\inf_{\|h-f\|_{\C^0}\leq \ve}\mathcal{H}^{d+p-m}\left(\mathcal{Z}_{W}^{h}\right)
	\eeq
	to be the smallest $(d+p-m)$-Hausdorff measure of $\mathcal{Z}_{W}^{h}$ among all functions 
	$h\in \mathcal{C}\left([0,1]^d,\R^m\right)$ with $\|h-f\|_{\C^0}\leq \ve$. In \cite{MN}, an upper bound on the number $\mathcal{N}^{f}_{W}(\ve)$  was recently established and  applied to provide quantitative estimates on the number of shock curves in  entropy weak solutions of scalar conservation laws with strictly convex fluxes.  Specifically, for $f\in\mathcal{C}^\alpha([0,1]^d,\R^m)$ with H\"older norm $\|f\|_{\mathcal{C}^{0,\alpha}}$ and $\ve>0$ sufficiently small, there exists a constant $C_W>0$ that depends only on $W$  such that 
	\bel{N-1}
	\mathcal{N}^f_{W}(\ve)~\leq~C_{W}\cdot  \left({\|f\|_{\mathcal{C}^{0,\alpha}}\over \ve}\right)^{m-p\over \alpha}.\eeq
	%where the constant $C_W>0$ depends only on $W$ and $\|f\|_{\mathcal{C}^{0,\alpha}}$ is the H\"older norm of $f$.
	% \begin{theorem}\label{Holder}
	%Assume that   $p+d\geq m$ and $f\in\mathcal{C}^\alpha([0,1]^d,\R^m)$ is H\"older continuous with  exponent $\alpha\in ]0,1]$. Then for every $\ve>0$ sufficiently small, it holds 
	%\[
	%\]
	%\bel{N-1}
	%\mathcal{N}^f_{W}(\ve)~\leq~C_{W}\cdot  \left({\|f\|_{\mathcal{C}^{0,\alpha}}\over \ve}\right)^{m-p\over \alpha},\eeq
	%where the constant $C_W>0$ depends only on $W$ and $\|f\|_{\mathcal{C}^{0,\alpha}}$ is the H\"older norm of $g$.
	%\end{theorem}
	The blow up rate $\left({1\over \ve}\right)^{m-p\over \alpha}$ with respect to $\ve$ is shown to be the best bound in  terms of power function  in \cite[Example 3.1]{MN} for a class of Lipscthiz functions $(\alpha=1)$ in the scalar case ($d=m=1$).  However, this still remains open for the multi-dimensional cases. Hence, the present paper aims to address the sharpness of  (\ref{N-1}) for general continuous function $f\in \C([0,1]^d,\R^m)$ with $d,m\geq 1$. In particular, we achieve the following lower quantitative estimate  for the class of H\"older continuous functions. 
	\begin{theorem}\label{M1} Assume that  $p< m\leq p+d$ and $W\subset\R^m$ is a $\C^1$-manifold  of dimension $p$. For every $0<\alpha\leq 1$ and $\lambda>0$, there exists a H\"older continuous function $f:[0,1]^d\to\R^m$ with  exponent $\alpha$ and the H\"older norm $\lambda$ such that
		\[
		\mathcal{N}^f_{W}(\ve)~\ge~C_{[W,\alpha,\lambda]}\cdot \left({1\over \ve\cdot 2^{4\cdot\sqrt{\alpha|\log_2\ve|}}}\right)^{m-p\over \alpha}	%\left({1\over 2^\sqrt{|\log_2\ve|}}\right)^{m-p\over \sqrt{\alpha}}
		\]
		for some constant $C_{[W,\alpha,\lambda]}>0$ that depends only  on $W$, $\alpha$, and $\lambda$.
	\end{theorem}
	Here the  constant $C_{[W,\alpha,\lambda]}$   is explicitly computed in Remark \ref{R}. Moreover, by using the concept of modulus of continuity and its inverse in Definition \ref{modu}, a  general result for continuous functions will be  proved  in Theorem \ref{Main} of Section \ref{2}. This  can  be easily extended  to the case of continuous functions $f:X \to Y$ where $X,Y$ are  smooth manifolds and $W\subseteq Y$ is a  smooth submanifold of $Y$. Finally, we remark that the factor $2^{4\cdot\sqrt{\alpha|\log_2\ve|}}$ in  Theorem \ref{M1} is necessary. Indeed,  we shall prove in the  Proposition \ref{Spn}  that  the  estimate  on $\mathcal{N}^{f}_{W}(\ve)$ in (\ref{N-1}) is not actually sharp for the case $\alpha=d=m=1$, $p=0$, and $W=\{0\}$. This leads to an open question on  the sharp estimate for $\mathcal{N}^f_{W}(\ve)$.

	\section{A lower bound on $\mathcal{N}^{f}_{W}(\ve)$}\label{2}
	In this section, we will establish a lower quantitative estimate  on the Hausdorff measure of $\mathcal{Z}^f_{W}$ for a constructed continuous $f\in\mathcal{C}\left([0,1]^d,\R^m\right)$ which admits  a given  modulus of continuity and    the set $W\subseteq \R^m$ being a  $\mathcal{C}^1$ manifold with $\mathrm{dim}(W)=p$. 
	For the sake of simplicity, we shall assume that $W$  consists of only one chart  $\R^m$, i.e., 
	\begin{itemize}
		\item [{\bf (A1)}.]  {\it There exists  a $\mathcal{C}^1$ diffeomorphism  $\phi$ between open subsets $U,V\subset\R^m$ such that $W\subset U$ and $\phi(W)=\R^p\times\{0\}\cap V$ and
			\bel{gamma}
			0~<~\gamma_{W}~\doteq~2\sqrt{m-p}\cdot\left({\sup_{x\in U}|\nabla\phi(x)|\over \inf_{x\in U}|\nabla\phi(x)|}\right)~<~\infty.
			\eeq	}
	\end{itemize} 
	For a general $\C^1$ manifold $W$ consists of multiple charts, one can just restrict the construction of  $f$ in a single chart of $W$ which has a smallest constant $\gamma_W$ among other charts. Toward to the main result, let us now recall some basic concepts on the modulus of continuity and its inverse.
	\begin{definition}\label{modu}  Given subsets $U\subseteq \R^d$ and $V\subseteq \R^m$, let $h:U\to V$ be continuous. The  minimal modulus of  continuity of $h$ is given by 
		\bel{omega}
		\omega_h(\delta)~=~\sup_{x,y\in U,|x-y|\leq \delta} |h(y)-h(x)|\qquad\forall \delta\in [0,\mathrm{diam}(U)].
		\eeq
		The inverse of the minimal modulus of  continuity of   $h$ is the map $s\to \Psi_h(s)$ is defined by 
		\bel{Psi}
		\Psi_{h}(s)~:=~\sup\left\{\delta\geq 0: |h(x)-h(y)|\leq s~~\forall |x-y|\leq \delta, x,y\in U\right\}
		\eeq
		for all $s\geq 0$.
	\end{definition}
	\noindent It is clear that $\Psi_h(s)=\infty$ for all $s\in [M_h,\infty[$ with $M_{h}:=\sup_{x,y\in U} |h(x)-h(y)|$.
	In particular, if $h$ is a constant function then $\Psi_{h}(s)=\infty$ for all $s\geq 0$. Otherwise, by the continuity of $h$, it holds
	\[
	\Psi_h(0)~=~0\qquad\mathrm{and}\qquad 0~<~\Psi_h(s)~\leq~\mathrm{diam}(U)\qquad\forall s\in ]0,M_h[.
	\]
	Moreover,  $\Psi_h(\cdot):[0,\infty[\to [0,\infty[$ is  increasing and superadditive
	\[
	\Psi_h(s_1+s_2)~\geq~\Psi_h(s_1)+\Psi_h(s_2)\qquad\forall s_1,s_2\geq 0.
	\]
	If the map $\delta\to \omega_h(\delta)$ is strictly increasing in $[0,\mathrm{diam}(U)[$ then $\Psi_h$ is the inverse of $\omega_{h}$, i.e., 
	\[
	\Psi_h(s)~=~\omega_{h}^{-1}(s)\qquad\forall s\in [0,M_h[.
	\]
	From the above observations, we  define a modulus of continuity as follows:
	\begin{definition}\label{beta} A function $\beta:[0,\infty]\to[0,\infty]$ is called a modulus of continuity if it is increasing, subadditive, and satisfies
		\[\lim_{\delta\to0+}\beta(\delta)~=~\beta(0)~=~0.\]
		We say that a continuous function $f:U\subset\R^d\to \R^m$ admits $\beta$ as a modulus of continuity if
		\begin{equation}
			\label{betadef}
			\sup_{x,y\in U,|x-y|\le s}|f(x)-f(y)|\le \beta(s)\qquad \forall s\ge0.
		\end{equation}
	\end{definition}	
	The main result in this paper is stated as follows: 	
	\begin{theorem}\label{Main} In addition to {\bf (A1)}, assume that  $p< m\leq p+d$.	For every modulus of continuity $\beta$, there exists a continuous function $f:[0,1]^d\to\R^m$ that admits $\beta$ as a modulus continuity and for $\ve>0$ sufficiently small
		\bel{lbman}
		\mathcal{N}^f_{W}(\ve)~\ge~\left(\frac{16}{\Psi_{\beta}(\gamma_{W}\ve)}\right)^{m-p}\cdot 2^{-4(m-p)\cdot\sqrt{\big|\log_2(\Psi_{\beta}(\gamma_W\ve))\big|}}.
		\eeq
		%for every $0<\ve<\frac{\beta(2^{-5})}{2\sqrt{m-p}}$. 
	\end{theorem}	
	{\bf Proof.} The proof is divided into three main  steps: 
	
	\underline{{\it Step 1}.} Consider the case $W=\{0\}$ and $p=0$. We claim that 
	\begin{itemize}
		\item [{\bf (G).}]  {\it There exists a continuous function $\tilde{f}:[0,1]^d\to\R^m$ that admits $\beta$ as a modulus of continuity and for every $0<\ve<{1\over2\sqrt{m}}\cdot\beta(2^{-5})$ it holds
			\bel{Lb-1}
			\mathcal{N}^{\tilde{f}}_{\{0\}}(\ve)~\geq~\left(\frac{16}{\Psi_{\beta}(2\sqrt{m}\ve)}\right)^m\cdot 2^{-4m\cdot\sqrt{\big|\log_2\big(\Psi_{\beta}\big(2\sqrt{m}\ve\big)\big)\big|}}
			\eeq
			with $\Psi_{\beta}$ being the inverse of the minimal modulus of  continuity of $\beta$.}
	\end{itemize}
	The construction of a desired function $\tilde{f}\in\mathcal{C}([0,1],\R^{m})$ in  {\bf (G).} will be done as follows: 
	
	{\bf 1.}  Let's first divide   $[0,1]$ into countably infinite  subintervals $[s_n,s_{n+1}]$ with
	\[
	s_1~=~0,\qquad s_n~=~\sum_{\ell=1}^{n}2^{-\ell}\qquad\forall n\geq 2.
	\]
	For every $n\geq 1$, we define  $u_n:[0,1]\to \R$ by
	\[
	u_n(s)~=~\sum_{k=0}^{2^{n^2}-1}c_n(s-s_n-4k\ell_n),\qquad \ell_n~=~2^{-n^2-n-2},
	\]
	where   $c_n:[0,1]\to\R$ is a sample function with $\mathrm{supp}(c_n)\subseteq [0,4\ell_n]$ such that  for all $s\in [0,2\ell_n]$
	\bel{cdef}
	c_n(s)~=~-c_n(4\ell_n-s)~=~{\beta(s)\over 2}\cdot \chi_{[0,\ell_n[}(s)+{\beta(2\ell_n-s)\over 2}\cdot \chi_{[\ell_n,2\ell_n]}(s)\,.
	\eeq
	%Consider the function $r:[0,1]\to\R$ such that 
	%\[
	%r(s)~=~\sum_{n=1}^{\infty}u_n(s)\qquad\forall s\in [0,1].
	%\]
	The   function $\tilde{f}=(\tilde{f}_1,\tilde{f}_2,\dots, \tilde{f}_{m})\in\mathcal{C}([0,1],\R^{m})$ is defined by 
	\bel{gi}
	\tilde{f}(x)~=~{1\over \sqrt{m}}\cdot \left(r(x_1),\dots, r(x_m)\right)\quad\forall x=(x_1,\dots, x_d)\in [0,1]^d
	\eeq
	with
	\[
	r(s)~\doteq~\sum_{n=1}^{\infty}u_n(s)\qquad\forall s\in [0,1].
	\]
	Since the modulus of continuity of $r$ is bounded by $\beta$, the modulus of continuity of $\tilde{f}$ is also bounded by $\beta$. Indeed, for every $s\geq 0$, one estimates
	\begin{eqnarray*}
		\omega_{\tilde{f}}(s)&=&\sup_{x,y\in [0,1]^d, |x-y|\leq s} |\tilde{f}(x)-\tilde{f}(y)|\\
		&=&\sup_{x,y\in [0,1]^d, |x-y|\leq s} {1\over \sqrt{m}}\cdot \left(\sum_{i=1}^{m}|r(x_i)-r(y_i)|^2\right)^{1\over 2}~\leq~\beta(s).
	\end{eqnarray*}
	Assume that  for every  $\ve>0$ satisfying
	\bel{con-ve}
	\ds{1\over2\sqrt{m}}\cdot\beta\left({\ell_{n_0+1}\over2}\right)~\le~\ve~\leq~ \ds {1\over2\sqrt{m}}\cdot\beta\left({\ell_{n_0}\over2}\right),
	\eeq 
	it holds 
	\bel{N-g-0}
	\mathcal{N}^{\tilde{f}}_{\{0\}}(\ve)~=~\inf_{\|g-f\|_{\mathcal{C}^0}\leq \ve}\mathcal{H}^{d-m}\left(\mathcal{Z}_{\{0\}}^{g}\right)~\geq~2^{mn_0^2}.
	\eeq
	In this case,  by the properties of an inverse of the minimal modulus of  continuity in (\ref{Psi}), we have that 
	\[
	\Psi_{\beta}(2\sqrt{m}\ve)~\geq~\Psi_{\beta}\left(\beta\left({\ell_{n_0+1}\over 2}\right)\right)~\geq~{\ell_{n_0+1}\over 2}~=~2^{-(n_0+1)^2-(n_0+1)-3}~\geq~2^{-(n_0+2)^2}.
	\] 
	Thus, one has
	\[
	n_0~\ge~ -2+\sqrt{-\log_2\Psi_{\beta}\big(2\sqrt{m}\ve\big)}
	\]
	and \eqref{Lb-1} follows from (\ref{N-g-0}).
	\medskip
	
	{\bf 2.} In the next two steps, we shall prove (\ref{N-g-0}). For every  $n\geq 1$ and  $k\in \{0,1,\dots, 2^{n^2}-1\}$, set
	\[
	a_{n,k}~=~s_n+(4k+1)\ell_n,\qquad b_{n,k}~=~s_n+(4k+3)\ell_n,
	\]
	we shall denote by 
	\bel{sq}
	\square_{n,\iota}~=~[a_{n,\iota_1},b_{n,\iota_1}]\times \cdots \times [a_{n,\iota_m},b_{n,\iota_m}]\quad\forall \iota\in \{0,1,\dots, 2^{n^2}-1\}^m.
	\eeq
	Fix $g\in \mathcal{C}\big([0,1]^d,\R^m\big)$ with   $\|\tilde{f}-g\|_{\C^0}\leq \ve$. By the definition of   $\mathcal{Z}_{\{0\}}^{g}$, we have 
	\[
	\mathcal{Z}_{\{0\}}^{g}~\supseteq~\bigcup_{n\geq 1, \iota\in \{0,1,\dots, 2^{n^2}-1\}^m}\left(\bigcup_{z\in [0,1]^{d-m}} \mathcal{Z}_{n,\iota}(z)\times \{z\}\right)
	\]
	with \[\mathcal{Z}_{n,\iota}(z)~=~\{y\in \square_{n,\iota}:g(y_1,...,y_m,z_1,...,z_{d-m})=0\}.\]
	Assume that for every  $1\leq n\leq n_0$ and $\iota\in \{0,1,\dots, 2^{n^2}-1\}^m$, the set 
	\bel{Con-k}
	\mathcal{Z}_{n,\iota}(z)~\neq~\varnothing\qquad\forall z\in [0,1]^{d-m}.
	\eeq
	In this case, we can bound the $(d-m)$-Hausdorff measure of $\mathcal{Z}_{\{0\}}^{g}$  by 
	\begin{eqnarray*}
		\mathcal{H}^{d-m}\left(\mathcal{Z}_{\{0\}}^{g}\right)&\ge&~\sum_{n=1}^{\infty}\sum_{\iota\in\{0,1,\dots, 2^{n^2}-1\}^m}\mathcal{H}^{d-m}\left(\bigcup_{z\in [0,1]^{d-m}} \mathcal{Z}_{n,\iota}(z)\times \{z\}\right)\\
		&\geq&~\sum_{n=1}^{n_0}\sum_{\iota\in\{0,1,\dots, 2^{n^2}-1\}^m}\mathcal{H}^{d-m}\left([0,1]^{d-m}\right)~=~\sum_{n=1}^{n_0}2^{m n^2}~\geq~2^{m n_0^2},
	\end{eqnarray*}
	and this yields  (\ref{N-g-0}).  
	
	\medskip

	{\bf 3.} To complete the proof, we need to verify (\ref{Con-k}). Fix $n\in \{1,\dots, n_0\}$, $\iota\in \{0,1,\dots, 2^{n^2}-1\}^m$, and $z\in [0,1]^{d-m}$, we consider the continuous map $h^z:\square_{n,\iota}\to\R^m$ such that 
	\bel{h}
	h^z(y)~=~y+\frac{\sqrt{m}\cdot \ell_n}{\beta(\ell_n/2)}\cdot g(y,z)\qquad\forall y\in \square_{n,\iota}.
	\eeq
	Notice that  $\square_{n,\iota}\subseteq [0,1]^m$ is a cube of size $2\ell_n$ centered at $c^{\iota,n}$ with 
	\[
	c^{\iota,n}_i~=~s_n+(4\iota_i+2)\ell_n\qquad\forall i\in\{1,2,\dots, m\}.
	\]
	Recall (\ref{gi}), (\ref{con-ve}), and $\|\tilde{f}-g\|_{\C^0}\leq \ve$, for every $y\in \square_{n,\iota}$ and $i\in\{1,2,\dots,m\}$, set $s:=y_i-s_n-4\iota_i\ell_n\in [\ell_n,3\ell_n]$, we estimate 
	\bel{invar1}
	\begin{split}
		\big|h^z_i(y)-c^{\iota,n}_i\big|~&=~\left|y_i+\frac{\sqrt{m}\cdot \ell_n}{\beta(\ell_n/2)}\cdot g_i(y,z)-s_n-(4\iota_i+2)\ell_n\right|\\
		~&\le~\frac{\sqrt{m}\cdot \ell_n}{\beta(\ell_n/2)}\cdot \ve+\left|y_i+\frac{\sqrt{m}\cdot \ell_n}{\beta(\ell_n/2)}\cdot f_i(y,z)-s_n-(4\iota_i+2)\ell_n\right|\\
		~&\le~{\ell_n\over 2}+\left|y_i+{\ell_n\over \beta({\ell_n/ 2})}r(y_i)-s_n-(4\iota_i+2)\ell_n\right|\\
		%~&=~\textcolor{blue}{{\ell_n\over 2}+\left|y_i+{\ell_n\over \beta({\ell_n/ 2)}}\cdot c_n(y_i-s_n-4\iota_i\ell_n)-s_n-(4\iota_i+2)\ell_n\right|.}\\%~\leq~{\ell}_n.}
		~&=~\frac{\ell_n}{2}+\left|s-2\ell_n+\ell_n\cdot\frac{c_n(s)}{\beta(\ell_n/2)}\right|\,.
	\end{split}
	\eeq
	By the definition of $c_n$ in (\ref{cdef}), both cases $s\in [\ell_n,2\ell_n]$ and $s\in [2\ell_n,3\ell_n]$ are similar, we shall bound  $\big|h^z_i(y)-c^{\iota,n}_i\big|$ for $s\in [\ell_n,2\ell_n]$. In this case, we have that
	\begin{equation*}
		\begin{split}
			\big|h^z_i(y)-c^{\iota,n}_i\big|~&=~\frac{\ell_n}{2}+\left|s-2\ell_n+\ell_n\cdot\frac{\beta(2\ell_n-s)}{2\beta(\ell_n/2)}\right|
		\end{split}
	\end{equation*}
	If $s\geq \ds {3\ell_n\over 2}$ then  $\big|h^z_i(y)-c^{\iota,n}_i\big|\leq \ds {\ell_n\over 2}+\ds \max\left\{2\ell_n-s,\ell_n\cdot\frac{\beta(2\ell_n-s)}{2\beta(\ell_n/2)}\right\}\leq \ell_n$. Otherwise, if $\ell_n\leq s< \ds {3\ell_n\over 2}$ then by the subadditivity of $\beta$, we have
	\[
	-\ell_n~=~-{\ell_n\over 2}- \left(\ell_n-\ell_n\cdot {\beta(\ell_n/2)\over 2\beta(\ell_n/2)}\right)~\leq~h^z_i(y)-c^{\iota,n}_i~\leq~{\ell_n\over 2}-{\ell_n\over 2}+\ell_n\cdot\frac{\beta(\ell_n)}{2\beta(\ell_n/2)}~\leq~\ell_n.
	\] 
	Thus, the map $y\mapsto h^z(y)$ is invariant in $\square_{n,\iota}$. Finally, by Brouwer's fixed point theorem, $h^z$ has a fixed point $y_z\in \square_{n,\iota}$, and (\ref{h}) implies that  $y_z$ belongs to the set $\mathcal{Z}_{n,\iota}(z)$ in (\ref{Con-k}). The proof of {\bf (G)}is complete.
	\v

	\underline{{\it Step 2.}} For every given $r_0>0$, we shall prove our result for the  case $W=[-r_0,r_0]^p\times\{0\}^{m-p}$. From {\bf (G)}, there exists a function  $\tilde{g}\in\mathcal{C}([0,1]^d,\R^{m-p})$  such that 
	\begin{itemize}
		\item $\tilde{g}$ admits $\beta$ as a modulus of continuity;
		\item For  every $0<\ve<\frac{1}{2\sqrt{m-p}}\cdot\beta(2^{-5})$, it holds
		\bel{lbmMM}
		\mathcal{N}^{g}_{\{0\}}(\ve)~\geq~\left(\frac{16}{\Psi_{\beta}(2\sqrt{m-p}\cdot\ve)}\right)^{m-p}\cdot 2^{-4(m-p)\cdot\sqrt{\big|\log_2\big(\Psi_{\beta}\big(2\sqrt{m-p}\cdot\ve\big)\big)\big|}}.
		\eeq
	\end{itemize}
	The continuous  function  $g:[0,1]^{d}\to \R^m$ defined by
	\[
	g(x)~=~(0,\tilde{g}(x))\qquad\forall x\in [0,1]^d,
	\]
	admits $\beta$ as a modulus of continuity. Moreover, if  $0<\ve\leq \min\left\{\frac{1}{2\sqrt{m-p}}\cdot\beta(2^{-5}),r_0\right\}$ then for every function $h=(h_1,...,h_m)\in \mathcal{C}([0,1]^d,\R^m)$  with $\|h-g\|_{\mathcal{C}^0}\le\ve$, it holds  
	$$
	h_i(x)~\in~[-r_0,r_0]\qquad\forall i\in \{1,\dots, p\}, x\in [0,1]^d.
	$$
	Thus, we  can bound the $(d-m+p)$ Hausdorff measure of $\mathcal{Z}_W^h$ by
	\begin{equation}
		\label{mMM}
		\begin{split}
			\mathcal{H}^{d-m+p}\left(\mathcal{Z}_W^h\right)~&=~\mathcal{H}^{d+m-p}\left(\left\{x\in[0,1]^d:h(x)\in[-r_0,r_0]^p\times\{0\}^{m-p}\right\}\right)\\
			~&=~\mathcal{H}^{d-m+p}\left(\{x\in[0,1]^d:(h_{p+1}(x),...,h_m(x))\in\{0\}^{m-p}\}\right)\\
			%~&=~\mathcal{H}^{d-m+p}\left(\mathcal{Z}_{\{0\}^{m-p}}^{(h_{p+1},...,h_m)}\right)\\
			~&\ge~\inf_{\|b-\tilde g\|_{\mathcal{C}^0}\le\ve}\mathcal{H}^{d-m+p}\left(\mathcal{Z}_{\{0\}}^b\right).
			%~&\ge~\left(\frac{16}{\Psi_{\beta}(2\sqrt{m-p}\cdot\ve)}\right)^{m-p}\cdot 2^{-4(m-p)\cdot\sqrt{\big|\log_2\big(\Psi_{\beta}\big(2\sqrt{m-p}\cdot\ve\big)\big)\big|}}
		\end{split}
	\end{equation}
	Substituting \eqref{lbmMM} into \eqref{mMM}, we obtain that
	\bel{lbmMM1}
	\begin{split}
		\mathcal{N}^{g}_{W}(\ve)~&=~\inf_{\|h-g\|_{\C^0}\le\ve}\mathcal{H}^{d-m+p}\left(\mathcal{Z}_W^h\right)~\geq~\inf_{\|b-\tilde{g}\|_{\mathcal{C}^0}\le\ve}\mathcal{H}^{d-m+p}\left(\mathcal{Z}_{\{0\}}^b\right)\\	
		~&\geq~\left(\frac{16}{\Psi_{\beta}(2\sqrt{m-p}\cdot\ve)}\right)^{m-p}\cdot 2^{-4(m-p)\cdot\sqrt{\big|\log_2\big(\Psi_{\beta}\big(2\sqrt{m-p}\cdot\ve\big)\big)\big|}}.
	\end{split}
	\eeq
	%	and this yields (\ref{lbman}) for $f=g$ and $W=[-r_0,r_0]^p\times\{0\}^{m-p}$. 
	
	\underline{{\it Step 3.}} To complete the proof,  we shall establish (\ref{lbman}) for a $\C^1$-smooth manifold  $W\subset \R^m$  satisfying {\bf (A1)}. Without loss of generality,  assume that for some $r_0>0$
	\[
	W_{r_0}~\doteq~[-\tilde{r}_0,\tilde{r}_0]^{p}\times\{0\}^{m-p} ~\subseteq~\phi(W),
	\]
	we consider $g$ for  $\ds r_0=\tilde{r}_0/ \lambda_2$ in Step 2 with 
	\bel{lamb}
	\ds\lambda_1~\doteq~\inf_{x\in U}|\nabla\phi(x)|\qquad\mathrm{and}\qquad \ds\lambda_2~\doteq~ \sup_{x\in U}|\nabla\phi(x)|.
	\eeq
	The desired function  $f:[0,1]^d\to\R^m$ is defined by
	\bel{ff}
	f(x)~=~\phi^{-1}\circ [\lambda_1\cdot g(x)]\qquad\forall x\in [0,1].
	\eeq
	Indeed, $f$ admits $\beta$ as a modulus of continuity since for every $x, y\in [0,1]^d$, it holds 
	\[
	|f(y)-f(x)|~\leq~{\lambda_1\cdot |g(y)-g(x)|\over \inf_{z\in U }|\nabla \phi (z)|}~=~ |g(y)-g(x)|.
	\]
	To verify \eqref{lbman}, let  $h\in\mathcal{C}([0,1]^d,\R^m)$ be such that $ \|h-f\|_{\C^0}\le\ve$. From (\ref{lamb}) and (\ref{ff}), one has that 
	\[
	\left\| {\phi\circ h\over \lambda_1}-g\right\|_{\C^0}~=~{1\over \lambda_1}\cdot \|\phi\circ h -\phi\circ f\|_{\C^0}~\leq~{\lambda_2\over \lambda_1}\cdot \|h-f\|_{\C^0}~\leq~{\lambda_2\ve\over \lambda_1},
	\]
	and this implies 
	\bel{est1}
	\begin{split}
		\mathcal{H}^{d-m+p}\left(\mathcal{Z}_W^h\right)&~=~\mathcal{H}^{d-m+p}\left(\left\{x\in [0,1]^d: (\phi\circ h)(x)\in \phi(W)\right\}\right)\\
		&~\geq~\mathcal{H}^{d-m+p}\left(\left\{x\in [0,1]^d: (\phi\circ h)(x)\in [-\tilde{r}_0,\tilde{r}_0]^{p}\times\{0\}^{m-p}\right\}\right)\\
		&~=~\mathcal{H}^{d-m+p}\left(\mathcal{Z}^{\phi\circ h}_{W_{r_0}}\right)~\geq~\mathcal{N}^{g}_{W_{r_0}}\left({\lambda_2\ve\over \lambda_1}\right).
	\end{split}
	\eeq
	Finally,  recalling (\ref{lbmMM1}), we get for every $h\in\mathcal{C}([0,1]^d,\R^m)$ with $ \|h-f\|_{\C^0}\le\ve$ that 
	\[
	\mathcal{H}^{d-m+p}\left(\mathcal{Z}_W^h\right)~\geq~\left(\frac{16}{\Psi_{\beta}(2\sqrt{m-p}\cdot \lambda_2\ve/\lambda_1)}\right)^{m-p}\cdot 2^{-4(m-p)\cdot\sqrt{\big|\log_2\big(\Psi_{\beta}\big(2\sqrt{m-p}\cdot \lambda_2\ve/\lambda_1\big)\big)\big|}},
	\]
	and (\ref{gamma}) yields (\ref{lbman}). 
	\endproof
	Notice that if $\beta(s)=\lambda s^{\alpha}$ for some $\lambda>0$ and $\alpha\in ]0,1]$ then from \eqref{Psi}, it holds
	\[
	\Psi_{\beta}(s)~=~\left(s\over \lambda\right)^{1\over \alpha}\qquad\forall s\in [0,\infty[.
	\] 
	In this case, we achieve an explicit estimate in (\ref{lbman}) by a direct computation. More precisely, we have the following remark.
	\begin{remark}\label{R}
		Under the same setting in Theorem \ref{Main}, if $\beta(s)=\lambda s^{\alpha}$ for some $\lambda>0$, $\alpha\in (0,1]$ then  there exists a H\"older continuous function $f:[0,1]^d\to\R^m$ with  exponent $\alpha$ and H\"older norm $\lambda$ such that 
		\[
		\mathcal{N}^f_{W}(\ve)~\ge~C_{[W,\alpha,\lambda]}\cdot \left(\frac{1}{\ve}\right)^{{m-p\over \alpha}}\cdot 2^{-{4(m-p)\over \alpha^{1/2}}\cdot\sqrt{\big|\log_2(\gamma_W\ve/\lambda)\big|}}.
		\]
		This particularly yields Theorem \ref{M1}.
	\end{remark}
	Finally,  we notice  that the factor $2^{-4m\cdot\sqrt{\big|\log_2\big(\Psi_{\beta}\big(2\sqrt{m}\ve\big)\big)\big|}}$ in  Theorem \ref{Main} is necessary. In other words,  the  estimate  on $\mathcal{N}^{f}_{W}(\ve)$ in (\ref{N-1}) is not actually sharp for the case $\alpha=d=m=1$, $p=0$, and $W=\{0\}$.
	
	\begin{proposition}\label{Spn}
		Assume that $d=m=1$, $p=0$, $W=\{0\}$ and  $\beta(s)=s$ for all $s\geq 0$. Then Theorem \ref{Main} does not hold if the factor $2^{-4(m-p)\cdot\sqrt{\big|\log_2(\Psi_{\beta}(\gamma_W\ve))\big|}}$ in \ref{lbman} is replaced by any positive constant.
	\end{proposition}
	%
	%\begin{remark} Assume that $d=m=1$, $p=0$, $W=\{0\}$ and  $\beta(s)=s$ for all $s\geq 0$. Then Theorem \ref{Main} does not hold if the factor $2^{-4(m-p)\cdot\sqrt{\big|\log_2(\Psi_{\beta}(\gamma_W\ve))\big|}}$ in \ref{lbman} is replaced by any positive constant.
	%\end{remark}
	{\bf Proof.} Arguing by contradiction, suppose that there exists a function $f\in\mathcal{C}([0,1],\R)$ and a constant $C_f\in (0,1]$ such that $f$ admits $\beta$ as a modulus of continuity and
	\begin{equation}
		\label{falsebd}
		\mathcal{N}_{\{0\}}^f~\ge~ \frac{C_f}{\ve}\qquad\forall \ve>0~\mathrm{small}.
	\end{equation}
	{\bf 1.} We first claim that 	for every $0<\ve\leq \ds {C_f^2\over 12}$ there exists $(y_i)_{i=1}^N\in [0,1]^N$ such that  
	\bel{N-y}
	N~\geq~\frac{C_f^2}{18\ve},\qquad |f(y_i)|~\geq~{\ve\over 2},\qquad |y_i-y_j|~\ge ~\frac{2\ve}{C_f}\quad \forall i\ne j.
	\eeq
	Indeed, dividing $[0,1]$ into $K_0=\lceil\frac{C_f}{3\ve}\rceil$ subintervals $[a_i,a_{i+1}]$ of length 
	\begin{equation}
		\label{len}
		\frac{2\ve}{C_f}~\le~ \ell_i~\le~\frac{3\ve}{C_f}\qquad\forall i\in \{0,\dots, K_0-1\},
	\end{equation}
	we consider a function $h_{\ve}\in\mathcal{C}([0,1],\R)$ which is defined in $[a_i,a_{i+1}]$ for every $i\in  \{0,\dots, K_0-1\}$ as follows:
	\begin{itemize}
		\item If $\ds \max_{x\in [a_i,a_{i+1}]}|f(x)|\le\frac{\ve}{2}$ then we set 
		\begin{equation}
			\label{hdef2}
			h_{\ve}(x)~=~\begin{cases}
				f(a_i)+x-a_i,&\ds a_i\le x\le a_i-f(a_i)+\frac{\ve}{2},\\[3mm]
				\ds \frac{\ve}{2},&\ds a_i-f(a_i)+\frac{\ve}{2}\le x\le a_{i+1}+f(a_{i+1})-\frac{\ve}{2},\\[3mm]
				f(a_{i+1})-x+a_{i+1},&\ds a_{i+1}+f(a_{i+1})-\frac{\ve}{2}\le x\le a_{i+1}.
			\end{cases}
		\end{equation}
		It is clear that  $h_{\ve}$ has at most $2$ zeros on $[a_i,a_{i+1}]$, and $\|h_{\ve}-f\|_{\C^0}\le\|h_{\ve}\|_{\C^0}+\|f\|_{\C^0}\le \ve$.
		\item Otherwise, if $\ds \max_{x\in [a_i,a_{i+1}]}|f(x)|>\frac{\ve}{2}$ then we divide $[a_i,a_{i+1}]$ into  $K_1=\lceil\frac{3}{C_f}\rceil$ subintervals $\big[a_i^j,a_i^{j+1}\big]$ of length at most $\ve$. For every  $j\in \{0,\dots, K_1-1\}$, we set 
		\begin{equation}
			\label{pwlin}
			h_{\ve}\big(\theta \cdot a_i^j+ (1-\theta)\cdot a_i^{j+1} \big)~=~\theta\cdot f\big(a_i^j\big)+\big(1-\theta)\cdot f(a_i^{j+1}\big),\qquad\theta\in [0,1].
		\end{equation}
		In this case, $h_{\ve}$ has at most $\frac{3}{C_f}+1$ zeros on $[a_i,a_{i+1}]$ and 
		\[
		\|h_{\ve}-f\|_{\C^0([a_i,a_{i+1}])}~\leq~\max_{0\leq j\leq K_1-1}\sup_{|x-y|\le a_i^{j+1}-a_i^j}|f(x)-f(y)|~\leq~\beta(\ve)~=~\ve.
		\]
	\end{itemize}
	Thus, set $\mathcal{I} =\big\{i\in \{0,\dots, K_0-1\}:\max_{x\in [a_i,a_{i+1}]}|f(x)|>\ve/2\big\}$ and $\eta=\#\I$. By  (\ref{N-f}) and (\ref{falsebd}), we have 
	\[
	\begin{split}
		{C_f\over \ve}&~\leq~\mathcal{H}^{0}\left(\mathcal{Z}_{\{0\}}^{h_\ve}\right)~\leq~\eta\cdot\left(\frac{3}{C_f}+1\right)+\left(K_0-\eta\right)\cdot 2~\leq~\eta\cdot\left(\frac{3}{C_f}+1\right)+\left({3\over C_f}+1-\eta\right)\cdot 2
	\end{split},
	\]
	and (\ref{falsebd}) yields 
	\[
	\eta~\geq~{C_f^2-6\ve-2\ve C_f\over (3-C_f)\ve}~\geq~{C_f^2\over 9\ve}.
	\]
	For every $i\in \I$, let $z_i\in\argmax_{x\in [a_i,a_{i+1}]} |f(x)|$ be such that $|f(z_i)|\geq \ve/2$. From the first inequality of (\ref{len}), one can pick a desired set of at least $N\geq\frac{C_f^2}{18\ve}$ points $y_i$ from the set $\{z_i:i\in\I \}$ which satisfies (\ref{N-y}).

	\medskip
	{\bf 2.} Using (\ref{N-y}), we  show that
	\begin{equation}
		\label{Nbd41}
		\mathcal{N}_{\{0\}}^f\left(\ve/4\right)~\le ~\left(1-\frac{7C_f^2}{36}\right)\cdot \frac{4}{\ve}\qquad\forall 0\leq \ve\leq \frac{C_f^2}{12}
	\end{equation}
	Divide $[0,1]$ into $K_{\ve}=\left\lceil{4\over \ve}\right\rceil$ subintervals $[b_k,b_{k+1}]$ with length smaller than $\ve/4$, let  $g_{\ve}:[0,1]\to\R$ be a continuous function such that  for all  $k\in \{0,\dots, K_{\ve}-1\}$, it holds
	\[
	g_{\ve}\big(\theta \cdot b_k+ (1-\theta)\cdot b_{k+1} \big)~=~\theta\cdot f\big(b_k\big)+\big(1-\theta)\cdot f(b_{k+1}\big),\qquad\theta\in [0,1].
	\]
	Up to a small variation, we can assume that $f(b_k)\neq 0$ for every $k\in \{1,\dots, K_{\ve}\}$ so that $g_{\ve}$ has at most one zero on each of the intervals $[b_k,b_{k+1}]$. By the construction, one has 
	\[
	\|g_{\ve}-f\|_{C^0}~\leq~\max_{0\leq k\leq K_{\ve}-1} |f(b_{k+1})-f(b_k)|~\leq~\max_{0\leq k\leq K_{\ve}-1} |b_{k+1}-b_k|~\leq~{\ve\over 4}.
	\]
	For every $k\in \{0,\dots, K_{\ve}-1\}$ such that $y_i\in [b_k,b_{k+1}]$ for some $i\in \{0,\dots, N\}$, it holds for all $x\in  [b_k,b_{k+1}]$ that 
	\[
	|g_{\ve}(x)|~\geq~|f(x)|-{\ve\over 4}~\geq~|f(y_i)|-|\beta(|x-y_i|)|-{\ve\over 4}~>~{\ve\over 4}- |b_{k+1}-b_k|~>~0.
	\]
	In this case,  $g_{\ve}$ is non-zero on the at least $N$ intervals $[b_{k}, b_{k+1}]$. Thus, we have 
	\[
	\mathcal{N}_{\{0\}}^f\left(\ve/4\right)~\leq~\mathcal{H}^0\left(\mathcal{Z}^{g_{\ve}}_{\{0\}}\right)~\leq~K_{\ve}-N~\leq~{4\over \ve}+1-{C_f^2\over 18\ve}
	\]
	and this yields (\ref{Nbd41}).
	\medskip
	
	{\bf 3.} Finally, applying (\ref{N-y}) $n$ times, we find that 
	\[
	\mathcal{N}_{\{0\}}^f\left({\ve\over 4^{n}}\right)~\leq~\left(1-\frac{7C_f^2}{36}\right)^n\cdot  {4^n\over \ve}.
	\]
	Thus, (\ref{falsebd}) does not holds for $\ve$ replaced by $\ds {\ve\over 4^n}$ with $n\geq 1$ sufficiently large so that $\left(1-\frac{7C_f^2}{36}\right)^n<C_f$. This concludes the proof.
	\endproof
	\v
	
	{\bf Acknowledgments.} This research by K. T. Nguyen was partially supported by National Science Foundation grant DMS-2154201.
	%
	%
	%		\begin{remark}
	%		Theorem \ref{main} can be extended non-smooth manifolds where the local chart is Lipschitz and has a Lipschitz inverse. In that case, $\lambda_1$ and $\lambda_2$ are replaced by Lipschitz constants of $\phi^{-1}$ and $\phi$ respectively.
	%	\end{remark}
	

\begin{thebibliography}{999}	
		\bibitem{AGN1} F. Ancona, O. Glass, and K.~T.~Nguyen,   Lower compactness estimates for scalar balance laws, {\it Comm. Pure Appl. Math}, 65(9), 1303-1329, 2012.
		
		\bibitem{AGN2} F. Ancona, O. Glass, and K.~T.~Nguyen,  On compactness estimates for hyperbolic systems of conservation laws, {\it Ann. Inst. H. Poincar\'e Anal. Non Lin\'eaire},  32(6), 1229-1257, 2015.
		
		\bibitem{AGN} F. Ancona, O. Glass, and K.~T.~Nguyen,  On Kolmogorov entropy compactness estimates for scalar conservation laws without uniform convexity, {\it SIAM J. Math. Anal.}, 51 (4), 3020--3051, 2019.
		\bibitem{ACN1} F. Ancona, P. Cannarsa, and K.~T.~Nguyen: Quantitative compactness for Hamilton Jacobi Equations, {\it Archive for Rational Mechanics and Analysis} 219 (2016), no. 2, 793--828
		\bibitem{ACN2}  F. Ancona, P. Cannarsa, and K.~T.~Nguyen: The compactness estimates for Hamilton Jacobi Equations depending on space, {\it Bulletin of the Institute of Mathematics, Academia Sinica} {\bf 11} (2016), 63--113
		\bibitem{BDN} S. Bianchini, P. Dutta, and K.~T.~Nguyen: Metric entropy for Hamilton-Jacobi equation with uniformly directional convex Hamiltonian,  {\it SIAM Journal on Mathematical Analysis} {\bf 54} (2022), 5551--5575.
		\bibitem{CDN} R. Capuani, P. Dutta Prerona, and K. T. Nguyen,  Metric entropy for functions of bounded total generalized variation, {\it SIAM J. Math. Anal.} 53 (1), 1168 -- 1190, 2021.
		
		\bibitem{MN}  A. Murdza and K.T. Nguyen, A quantitative version of the transversality theorem, Communications in Mathematical Sciences, to appear.
		\bibitem{Th1} J.M. Bloom, The local structure of  maps of manifolds. B.A. thesis, Harvard 2004. (Link: https://www.math.harvard.edu/media/ThesisXFinal.pdf)
		
		\bibitem{Thom} Ren\'e Thom, Quelques propri\'et\'es globales des vari\'et\'es differentiables, {\it Commentarii Mathematici Helvetici} {\bf 28} (1954), no. 1, 17-86.
		\bibitem{N1} R. M.~Hardt and T Rivi\`ere, Connecting rational homotopy type singularities, {\it Acta Math.} {\bf 200} (2008), 15-83.
		\bibitem{N2} S.~Li, A note on generic transversality of Euclidean submanifolds, {\it Manuscripta Math.} {\bf 162} (2020), 213-219.
	\end{thebibliography}
\end{document}